\documentclass{article}
\usepackage{amsmath}
\usepackage{amssymb}
\usepackage{amsfonts}

  \usepackage{paralist}
  \usepackage{graphics} %% add this and next lines if pictures should be in esp format
  \usepackage{epsfig} %For pictures: screened artwork should be set up with an 85 or 100 line screen
 \usepackage[colorlinks=true]{hyperref}
\hypersetup{urlcolor=blue, citecolor=red}

\newtheorem{theorem}{Theorem}[section]

\newtheorem{definition}[theorem]{Definition}

\newcommand{\rf}[1]{(\ref{#1})}
\newcommand{\ods}{\par \vspace{1em} \par}

\newcommand{\ba}{\begin{array}}
\newcommand{\ea}{\end{array}}

\newcommand{\be}{\begin{equation}}
\newcommand{\ee}{\end{equation}}

\newcommand{\const}{{\rm const}}

\newcommand{\R}{{\mathbb R}}  
\newcommand{\C}{{\mathbb C}}  
\newcommand{\T}{{\mathbb T}}  
\newcommand{\Z}{{\mathbb Z}}

\newcommand{\av}[1]{\mbox{$\langle #1 \rangle$}}

\newcommand{\dis}{\displaystyle }

\newenvironment{Proof}{\par \vspace{2ex} \par
\noindent \small {\it Proof:}}{\hfill $\Box$ 
\vspace{2ex} \par }

\newcommand{\no}{ \noindent }

\begin{document}

\title{ \scshape Some implications of a new definition of the exponential function on time scales}

\author{
 { \scshape Jan L.\ Cie\'sli\'nski}\thanks{\footnotesize
 e-mail: \tt janek\,@\,alpha.uwb.edu.pl} 
\\ {\footnotesize Uniwersytet w Bia{\l}ymstoku,
Wydzia{\l} Fizyki, ul.\ Lipowa 41, 15-424
Bia{\l}ystok, Poland}
}

\date{}

\maketitle

\bigskip

\begin{abstract}
We present a new approach to exponential functions on time scales and to timescale analogues of ordinary differential equations. We describe in detail the Cayley-exponential function and associated trigonometric and hyperbolic functions. We show that the Cayley-exponential is related to implicit midpoint and trapezoidal rules, similarly as  delta and nabla exponential functions are related to Euler numerical schemes. Extending these results on any Pad\'e approximants, we obtain Pad\'e-exponential functions.  Moreover, the exact exponential function on time scales is defined. 
Finally, we present applications of the Cayley-exponential function in the $q$-calculus and suggest a general  approach to dynamic systems on Lie groups. 
\end{abstract}

\ods



\ods

\section{Introduction} 
Time scales were introduced in order to unify differential and difference calculus \cite{Hi,Hi2}. A time scale $\T$ is any non-empty closed subset of $\R$.  In this paper we consider the problem of finding good analogues of exponential, hyperbolic and trigonometric functions on time scales, and also finding $\T$-analogues of ordinary differential equations (ODEs). 
In the first part of this paper we present results obtained in \cite{Ci-trig}. Then, many new developments are briefly described. In particular, we show how to generate  $\T$-analogues of ODEs using known numerical schemes, including the implicit midpoint rule, the trapezoidal rule, the discrete gradient method and Pad\'e approximants. We introduce exact $\T$-analogues of elementary and special functions and exact $\T$-analogues of ODEs. We also present some new results in the $q$-calculus,  see \cite{Ci-q}. Finally, we show how to define $\T$-analogues of functions on Lie groups using the Cayley transformation. 

The classical 
Cayley transformation, \ $\dis  z \rightarrow {\rm cay} (z,a) = (1 - a z)^{-1} (1 + a z)$, \  is a conformal transformation of the complex plane. Generalizations of the Cayley transformation on Lie groups and operator spaces are well known, see, e.g., \cite{Is-Cay}. 

The Cayley transformation maps the imaginary axis into the unit circle. The Cayley-exponential function and other new exponential functions presented in this paper also have this property. Therefore related $\T$-analogues of trigonometric functions satisfy exactly the Pythagorean identity: $\cos^2 t + \sin^2 t = 1$.

\section{Preliminaries and notation} In this section we briefly present some preliminaries on time scale calculus, with a special stress on standard approaches to exponential, hyperbolic and trigonometric functions. 

 \subsection{Basic notation and definitions}

\begin{itemize}
\item Forward jump operator $\sigma$: \quad \ 
$
\sigma (t) = \inf \{ s \in \T : s > t \} \equiv t^\sigma$
 \item  Backward jump operator $\rho$: \quad $  \rho (t) = \sup \{ s \in \T : s < t \}$  
\item 
Right-dense points:   \    $\sigma (t) = t$.  \quad Left-dense points:   \   $\rho (t) = t$.
\item
Right-scattered points:   \   $\sigma (t) > t$.  \quad Left-scattered points:   \   $\rho (t) < t$.
\item 
Graininess: \quad  $\dis  \mu (t) = \sigma (t) - t $. 
\item 
\emph{Rd-continuous} function is, by definition, continuous at right-dense points and has a finite limit at left-dense points. 
\item Graininess $\mu$ is not continuous at points which are left-dense and right-scattered but  
 is always rd-continuous.

\item There are two main $\T$-analogues of the $t$-derivative: 
\[  \ba{ll} 
\emph{Delta derivative} & \qquad \dis
f^\Delta (t ) : = \lim_{ \ba{c} s \rightarrow t \\ s \neq \sigma (t) \ea} \frac{ f ( \sigma (t) ) - f ( s )}{\sigma (t) - s}  \ , \\[4ex] 
\emph{Nabla derivative} & \qquad \dis
 f^\nabla (t) : = \lim_{ \ba{c} s \rightarrow t \\ s \neq \rho (t) \ea} \frac{ f ( \rho (t) ) - f ( s )}{\rho (t) - s} \ .
\ea \]

\end{itemize}

\subsection{Exponential functions}
 
 \emph{ Delta exponential function} (see \cite{Hi}), denoted by  $e_\alpha (t, t_0)$, is the unique solution of the  initial value problem 
\be  \label{ini-del}
\quad x^\Delta = \alpha (t) \, x \ , \qquad x (t_0) = 1 \ ,
\ee
where $\alpha: \T \rightarrow \C$  is a given function. 
\emph{ Nabla exponential function} (see \cite{ABEPT}),  denoted by \  ${\hat e}_\alpha (t, t_0)$, satisfies: 
\[
\quad x^\nabla = \alpha (t) \, x \ , \qquad x (t_0) = 1 \ .
\]
In the continuous case both exponential functions reduces to 
\be  \label{cont-exp}
\T = \R \quad  \Longrightarrow \quad  e_\alpha (t, t_0) = {\hat e}_\alpha (t, t_0) = \exp \int_{t_0}^t \alpha (s) \, d s \ .
\ee
In particular, 
\[
\T = \R \ , \ \alpha (t) = z  \quad \Longrightarrow \quad   \dis e_\alpha (t) = {\hat e}_\alpha (t) = e^{z t} \ . 
\]
In the discrete  constant case ($\T = h \Z$, \ $\alpha (t) = z \in \C$) we have 
\[
 e_z (t) = \left(1 + \frac{z t}{n} \right)^n  , \qquad  {\hat e}_z (t) = \left(1 - \frac{z t}{n} \right)^{-n}  , \qquad  t = n h  \ .
\]
Another exponential function is related to the so called diamond-alpha derivative \cite{MT}.  However, the associated differential equation is of second order.

\subsection{ Hyperbolic and trigonometric functions on time scales} 
\label{sec-hyptrig}

There exist two different approaches to hyperbolic and trigonometric functions on time scales. This ambiguity can be explained as follows. In the continuous case we have
\[
\T = \R \quad \Longrightarrow \quad \cos x = \frac{e^{i x} + e^{- i x}}{2}  \ , \quad  \sin x = \frac{e^{i x} - e^{- i x}}{2 i}  \ .
\]
Unfortunately, the property $e^{- i x} = (e^{i x})^{-1}$, crucial for some important properties of trigonometric functions, does not extend on delta and nabla exponential functions: 
\[
e_{-\alpha} (t, {t_0}) \neq e_\alpha^{-1} (t, {t_0}) \ , \qquad 
{\hat e}_{-\alpha} (t, {t_0}) \neq {\hat e}_\alpha^{-1} (t, {t_0}) \ .
\]
Therefore, trying to extend definitions of trigonometric functions on time scales, one has to choose between two natural possibilities: to replace $e^{-i x}$ either by $(e_{ix})^{-1}$ or by $e_{- i x} $. 
Hilger chose the  first option \cite{Hi-spec}, i.e., 
\[
 \cosh_\alpha (t) = \frac{ e_\alpha (t) + e_\alpha^{-1} (t) }{2}  \ , \qquad
 \sinh_\alpha (t) = \frac{ e_\alpha (t) - e_\alpha^{-1} (t) }{2} \ .
\]
The main advantage of this definition is the identity
\[
 \cosh_\alpha^2 (t) - \sinh_\alpha^2 (t) = 1  \ .  
\] 
Unfortunately, derivatives \ $\cosh_\alpha^\Delta (t)$, $\sinh_\alpha^\Delta (t)$,   \ are not proportional to \  $\sinh_\alpha (t)$ and $\cosh_\alpha (t)$, respectively. Instead, we have
\[ \ba{l} \dis
\cosh_\alpha^\Delta (t) = \frac{\frac{1}{2} \mu \alpha^2}{1 + \mu \alpha} \cosh_\alpha (t) + \frac{\alpha + \frac{1}{2} \mu \alpha^2 }{1 + \mu \alpha} \sinh_\alpha (t) \ , \\[3ex] \dis
\sinh_\alpha^\Delta (t) = \frac{\frac{1}{2} \mu \alpha^2}{1 + \mu \alpha} \sinh_\alpha (t) + \frac{\alpha + \frac{1}{2} \mu \alpha^2 }{1 + \mu \alpha} \cosh_\alpha (t) \ ,
\ea \]
where $\mu, \alpha$ may depend, in general, on $t$. Another disadvantage is \ $\cosh_{i \omega} (t) \notin \R$ and  $\sinh_{i \omega} (t) \notin i \R$  \  (for $\omega\in \R$).  Therefore, it is difficult to extend this approach on trigonomeric functions. However, after some algebraic considerations Hilger arrived at surprisingly simple formulae (see \cite{Hi-spec})
\be
\omega (t) = \const \qquad \Longrightarrow \qquad \cos_\omega (t) := \cos (\omega t) \ , \quad \sin_\omega (t) := \sin (\omega t) \ .
\ee
In section~\ref{sec-exact} the above definition of trigonometric functions on time scales will be interpreted in a broader framework of \emph{exact discretizations}.

The second (more popular) approach to hyperbolic and trigonometric functions has been proposed by Bohner and Petersson \cite{BP12,BP-I}
  \be   \ba{l} \dis \label{trig-BP}
 \cosh_\alpha (t) = \frac{e_\alpha (t) + e_{-\alpha} (t) }{2} \ ,  \qquad \cos_\omega (t) = \cosh_{i\omega} (t) \ , \\[3ex] \dis
  \sinh_\alpha (t) = \frac{e_\alpha  (t) - e_{-\alpha} (t) }{2} \ ,  \qquad  \sin_\omega (t) = - i \sinh_{i\omega} (t) \ . 
\ea \ee
We see that \ $\cos_\omega (t) \in \R$ \ and  \ $\sin_\omega (t) \in \R$ \ for  
 $\dis \omega (t) \in \R$. What is more, the derivatives of hyperbolic and trigonometric functions are good analogues of the continuous case:  
\[  \ba{l}
\sinh_\alpha^\Delta (t) = \alpha \cosh_\alpha (t) \ , \quad   \sin_\omega^\Delta (t) = \omega \cos_\omega (t) \ ,  \\[3ex]
\cosh_\alpha^\Delta (t) = \alpha \sinh_\alpha (t) \ ,   \quad  \cos_\omega^\Delta (t) = - \omega \sin_\omega (t) \ . 
\ea \]
Unfortunately, in  place of Pytha\-go\-rean identities we have qualitatively different equalities, namely 
\[ 
\cosh_\alpha^2 (t) - \sinh_\alpha^2 (t) = e_{-\mu\alpha^2} (t)  , \qquad 
\cos_\omega^2 (t) + \sin_\omega^2 (t) = e_{\mu\omega^2} (t) \ .
\]
Therefore, functions $\cos_\omega (t)$ and $\sin_\omega (t)$, defined by \rf{trig-BP}, are not bounded.

\section{ The Cayley-exponential function on $\T$}
\label{sec-Cexp}

We are going to define the Cayley-exponential function as a solution of an appropriate initial value problem. The definition presented in \cite{Ci-trig} is equivalent, see  Theorem~\ref{cyl}.

\subsection{ Definition of the C-exponential function}

\begin{definition} \label{def-Cay}
The Cayley-exponential (C-exponential) function  $E_\alpha (t, t_0)$ satisfies the following initial value problem:
\be  \label{ini-Cay}
       x^\Delta (t) = \alpha (t) \, \av{x (t)}  \ ,  \qquad  x (t_0) = 1 \ , 
\ee
\no where $\alpha$ is regressive (i.e., $\mu \alpha \neq \pm 2$) and rd-continuous on $\T$, and \ 
\[
 \av{x (t)}  = \frac{ x (t) + x (t^\sigma) }{2} \ , \qquad {\rm shortly:} \quad \av{x} = \frac{x + x^\sigma}{2} \ . 
\]
Moreover, we denote $E_\alpha (t) := E_\alpha (t, 0)$. 

\end{definition}

\no In the continuous case the C-exponential function reduces to the usual exponential, see \rf{cont-exp}. In the discrete case we have
\[
\T = h \Z , \  \alpha = \const , \quad  \Longrightarrow \quad  
  E_\alpha (t) = \left( \frac{1 + \frac{1}{2 n} t \alpha }{1 - \frac{1}{2 n} t \alpha }  \right)^n \ , \  t = n h  \ .
\]
Similar formulas in the \emph{discrete} case were known a long time ago to Ferrand \cite{Fer} and Duffin \cite{Duf}, and then appeared several times in different contexts \cite{BMS,DJM1,Is-Cay,NQC,ZD}. They are clearly related to the Cayley transformation, see   \cite{Is-Cay}. 

\subsection{ Properties of the C-exponential function}

\begin{theorem} \label{bij}
 There is a bijection between Cayley exponential functions and delta exponential functions. Namely,  \  $E_\alpha (t, t_0) = e_\beta (t, t_0)$, \ if    
\be  \label{betal}
 \alpha (t) = \frac{\beta (t)}{1 + \frac{1}{2} \mu (t) \beta (t)} \ , \qquad
\beta (t) = \frac{ \alpha (t)}{1 - \frac{1}{2}{\mu (t) \alpha (t)}} \ ,
\ee  
where  $\mu \alpha \neq \pm 2$ and $\mu \beta \neq - 1$. 
\end{theorem}

\begin{Proof}
Suppose that $x = x (t)$ satisfies \rf{ini-Cay}. Using \ $x^\sigma = x + \mu x^\Delta$ \  we obtain
\[
 (1 - \frac{1}{2} \alpha \mu ) \ x^\Delta = \alpha x \ .
\]
Hence, $x (t) = e_\beta (t, t_0)$, where $\beta$ is given by \rf{betal}. The assumption $\mu\alpha \neq \pm 2$ (see Definition~\ref{def-Cay}) quarantees that $\beta$ exists and $\mu \beta + 1 \neq 0$ ($\mu$-regressivity, see \cite{Ci-trig}). 

Likewise, suppose that $x = x (t)$ satisfies $x^\Delta = \beta x$, $x (t_0) = 1$. Substituting  \ $x = x^\sigma - \mu x^\Delta$ \  we get $(1 + \mu \beta) x^\Delta = \beta x^\sigma$. Adding both equations we obtain 
\[
    ( 2 + \mu \beta ) x^\Delta = 2 \beta \av{x} \ .
\]
Therefore, $x (t) = e_\alpha (t, t_0)$, where $\alpha$ is given by \rf{betal}. The assumption $\mu\beta \neq - 1$ quarantees that $\mu \alpha \neq - 2$.  
\end{Proof}

\begin{theorem}  \label{cyl}
The Cayley-exponential function $E_\alpha$ can be expressed as follows 
\be  \label{Ealpha}
  E_\alpha (t, t_0) := \exp \left( \int_{t_0}^t \zeta_{\mu (s)} ( \alpha (s)) \Delta s \right) \ ,  
\ee
where 
$\dis 
  \zeta_\mu (z) :=  \frac{1}{\mu} \log \frac{1 + \frac{1}{2} z \mu}{1 - \frac{1}{2} z \mu}$  \  (for $\mu \neq 0$),  i.e.,  $\dis  z = \frac{2}{\mu} \tanh \frac{\mu \zeta}{2} $,   \ and  \ $\zeta_0 (z) := z$.
 
\end{theorem}

\begin{Proof}
Using Theorem~\ref{bij} we can apply Hilger's results (the cylinder transformation) \cite{Hi,Hi-spec} and the formula \rf{Ealpha} follows in a straightforward way.  
\end{Proof}

The following properties of the Cayley-exponential function have been derived in \cite{Ci-trig}. The first formula is clearly related to the Cayley transformation. 
\[  \ba{l}  \dis
\displaystyle E_\alpha (t^\sigma, t_0) = \frac{1 + \frac{1}{2} \mu (t) \alpha (t)}{1 - \frac{1}{2} \mu (t) \alpha (t) } \ E_\alpha (t, t_0) \ , \qquad 
E_\alpha (t, t_0) \, E_\alpha (t_0, t_1) = E_\alpha (t, t_1) \ , \\[4ex] \dis
\overline{ E_\alpha (t, t_0)} = E_{\bar \alpha} (t, t_0) \ , \qquad 
( E_\alpha (t, t_0) )^{-1} = E_{-\alpha} (t, t_0) \ , \\[3ex] \dis
E_\alpha (t, t_0) \ E_\beta  (t, t_0) =  E_{\alpha\oplus \beta} (t, t_0) \ , \quad {\rm where} \qquad 
 \alpha \oplus \beta := \frac{\alpha + \beta}{1 + \frac{1}{4} \mu^2  \alpha \beta }  .
\ea \] 
Interestingly enough, the last formula is identical with the Lorentz velocity transformation of special relativity  provided that we interpret $\frac{2}{\mu}$ as the speed of light.

\subsection{Numerical advantages of the C-exponential function } 

In the continuous case $e_\alpha, {\hat e}_\alpha$ and $E_\alpha$ become identical, see \rf{cont-exp}. We are going to compare their ``accuracy'' at a right-scattered point $t$ 
(i.e., $t^\sigma - t = \mu \neq 0$), assuming $\alpha (t^\sigma) = \alpha (t) = \alpha$. 
\[  \ba{l} \dis
   E_\alpha (t^\sigma, t) = \frac{1 + \frac{1}{2} \mu \alpha }{1 - \frac{1}{2} \mu  \alpha  } = 1 + \alpha \mu + \frac{1}{2} (\alpha \mu)^2 + \frac{1}{4} (\alpha \mu)^3  + \ldots \ ,  
\\[3ex] \dis
e_{\alpha} (t^\sigma, t) = 1 + \alpha \mu \ , \qquad 
{\hat e}_\alpha (t^\sigma, t) = \frac{1}{1 - \alpha \mu} = 1 + \alpha \mu +  (\alpha \mu)^2 + \ldots 
\ea \]
where in the second formula $\alpha$ is evaluated at $t^\sigma$. It means that our assumpation ($\alpha^\sigma = \alpha$) is essential. Comparing the above expansions with the continuous case
\[
\quad \exp (\alpha \mu) = 1 + \alpha \mu + \frac{1}{2} (\alpha \mu)^2 + \frac{1}{6} (\alpha \mu)^3 +\ldots \ ,
\]
we conclude that \ $E_\alpha (t^\sigma, t)$  is a second-order approximation of $\exp (\alpha \mu)$, while $e_\alpha (t^\sigma, t)$ and ${\hat e}_\alpha (t^\sigma, t)$ are  approximations of the first order.

\section{Cayley-hyperbolic and Cayley-trigonometric functions on $\T$}

A direct consequence of Definition~\ref{def-Cay} is  better theory of hyperbolic and trigonometric fuctions. Cayley-hyperbolic and Cayley-trigonometric functions are defined in a natural way: 
\be \ba{l} \dis
{\rm Cosh}_\alpha (t) := \frac{E_\alpha (t) + E_{-\alpha} (t) }{2} \ , \quad {\rm Sinh}_\alpha (t) := \frac{E_\alpha (t) - E_{-\alpha} (t) }{2} \ , 
\\[3ex] \dis 
{\rm Cos}_\omega (t) := \frac{E_{i \omega} (t) + E_{- i \omega} (t) }{2} \ , \quad 
{\rm Sin}_\omega (t) := \frac{E_{i \omega} (t) - E_{- i \omega} (t) }{2 i} \ .
\ea \ee
We have no other choice because the Cayley-exponential function enjoys good properties like \ $(E_\alpha (t) )^{-1} = E_{-\alpha} (t)$,   $E_{\bar \alpha} (t) = \overline{ E_\alpha (t)  }$ and 
$|E_{i \omega } (t)| = 1$ for $\omega (t) \in \R$. 

 C-hyperbolic and C-trigonometric functions  combine advantages of both  Hilger's and Bohner-Peterson's approach.  

\begin{theorem} We assume that $\alpha$, $\omega$ are rd-continuous and $\alpha \mu \neq \pm 2$, $\omega \mu \neq \pm 2 i$. Then 
\[ \ba{ll} \dis
  {\rm Cosh}_\alpha^2 (t) - {\rm Sinh}_\alpha^2 (t) = 1 \ , & \quad {\rm Cos}_\omega^2 (t) + {\rm Sin}_\omega^2 (t) = 1 \ , 
\\[3ex]\dis 
{\rm Cosh}_\alpha^\Delta (t) = \alpha (t) \ \av{ {\rm Sinh}_\alpha (t)} \ , & \quad
{\rm Sinh}_\alpha^\Delta (t) = \alpha (t) \ \av{ {\rm Cosh}_\alpha (t)} \ , 
\\[3ex]\dis
{\rm Cos}_\omega^\Delta (t) = - \omega (t) \ \av{ {\rm Sin}_\omega (t)} \ , & \quad
{\rm Sin}_\omega^\Delta (t) = \omega (t) \ \av{ {\rm Cos}_\omega (t)} \ . 
\ea \]
\end{theorem}

\begin{Proof}  Straightforward calculation (see also \cite{Ci-trig}).
\end{Proof}

\section{ $\T$-analogues of ODE  motivated by numerical sche\-mes} 
\  Constructing $\T$-analogues of ordinary differential equations (ODEs) one usually replaces $t$-deri\-vatives by delta derivatives \cite{ABOP,BP-I} or, less frequently, by nabla derivatives \cite{ABEPT}. Thus \ $x^\Delta = \alpha x$ \ and \ $x^\nabla = \alpha x$ \ are standard $\T$-analogues of  \ $\dot x = \alpha x$. 
The results of section~\ref{sec-Cexp} suggest, however, that equation \ $x^\Delta = \alpha \av{x}$ \ is another (perhaps even better) $\T$-analogue of this equation. Certainly $\T$-analogues are not unique. 

In this section we are going to show a correspondence between numerical schemes and $\T$-analogues of a given ODE.  
We consider a general ODE:
\[
 \dot x = f (x, t) \ , \qquad t \in \T \ , \qquad x (t) \in \C^N \ , \qquad f (x (t), t) \in \C^N \ ,
\] 
and present its $\T$-analogues corresponding to several numerical methods. 

\subsection{Euler schemes} 

The most popular $\T$-analogues are associated with Euler methods. 
Delta dynamic equations correspond to the forward (explicit) Euler scheme, while nabla dynamic equations corespond to backward (implicit) Euler scheme
\[ \ba{ll} 
\emph{ Forward\ Euler\ scheme } & \qquad x^\Delta (t) = f (x (t), t) \ , \\[2ex]
\emph{ Backward Euler scheme}  & \qquad  x^\nabla (t) = f (x (t), t) \ .
\ea \]

\subsection{  Trapezoidal rule}  In the autonomous case we have
\[ \ba{ll}
\emph{Trapezoidal rule} & \qquad x^\Delta = \frac{1}{2} \left( f (x) + f (x^\sigma) \right) \ .
\ea \]
In the non-autonomous case we can consider at least two different possibilities:
\[ \ba{l}
x^\Delta = \frac{1}{2} \left( f (x, t) +  f (x^\sigma, t^\sigma) \right) \ , \\[2ex]
x^\Delta = \frac{1}{4} \left( f (x, t) + f (x^\sigma, t) + f (x, t^\sigma) + f (x^\sigma, t^\sigma) \right) 
\ea \] 
The first one is related to the standard trapezoidal scheme, the second one is more symmetric and has some advantages. Taking \ $f (x, t) = \alpha (t)  x$ \ we have, respectively, 
\[
  x^\Delta = \av{\alpha x} \    \dis \Rightarrow \    x^\sigma = \frac{1 + \frac{1}{2} \mu \alpha}{1 - \frac{1}{2} \mu \alpha^\sigma} \ x \ ,  \qquad
x^\Delta = \av{\alpha}  \av{x} \   \dis  \Rightarrow \     x^\sigma = \frac{1 + \frac{1}{2} \mu \av{\alpha} }{1 - \frac{1}{2} \mu \av{ \alpha } }  \ x \ . 
\]
These equations define next new exponential functions on time scales. The second of these functions  has better properties because it maps the imaginary axis into the unit circle for any function $\alpha = \alpha (t)$.

\subsection{  Implicit midpoint rule}  In the autonomous case we have 
\[ \ba{ll}
\emph{Implicit midpoint rule} & \qquad \dis x^\Delta =  f \left( \frac{x + x^\sigma}{2}\right) \ .
\ea \]
It is not clear how to extend this formula on the non-autonomous case, because in the general case  $\frac{1}{2} (t + t^\sigma)\notin \T$.  One may consider, for instance, the following scheme
\[
 x^\Delta = \frac{1}{2} \left( f \left( \frac{x + x^\sigma}{2}, t  \right) + f \left( \frac{x + x^\sigma}{2}, t^\sigma  \right)  \right) \ . 
\]

\subsection{ Discrete gradient method } 
This is an energy preserving numerical scheme for Hamiltonian systems \cite{LaG,MQR2}. Here we consider   one-dimensional separable systems  \ $H (p, q) = T (p) + V (q)$. The equations of motion read
\[
 \dot q = \frac{\partial T}{\partial p} \ , \ \qquad  \dot p = - \frac{\partial V}{\partial q} \ .
\]
The discrete gradient method yields the following $\T$-analogue of these equations of motion:
\be  \label{disgrad}
 q^\Delta = \frac{\Delta T}{\Delta p} \ , \qquad p^\Delta = - \frac{\Delta V}{\Delta q} \ .
\ee
where the ``discrete gradient'' is defined as \par
\[  \dis
\frac{\Delta T}{\Delta p} (p) := \lim_{P  \rightarrow p} \frac{ T(p^\sigma) - T(P) }{p^\sigma - P} \ , \qquad 
\ \frac{\Delta V}{\Delta q} (q) := \lim_{Q  \rightarrow q} \frac{ V(q^\sigma) - V(Q) }{q^\sigma - Q} \ .
\]
The dynamic system \rf{disgrad} preserves exactly the energy, i.e., $H^\sigma = H$ (more detailed discussion will be presented elsewhere). We point out that equations \rf{disgrad} differ from equations of motion considered in \cite{ABR}.

\subsection{Classical harmonic oscillator equation}  
Implicit midpoint, trapezoidal and discrete gradient schemes yield the same  $\T$-analogue of the harmonic oscillator equation ($\ddot q + \omega_0^2 q = 0$):
\[
 q^{\Delta\Delta} + \omega_0^2 \ \av{\av{q}} = 0 \ , \qquad  
\av{\av{q}}  := \frac{ q^{\sigma\sigma} + 2 q^\sigma + q }{4} \ . 
\] 
Of course, all these schemes yield different results for other, nonlinear, equations.

If \ $\omega (t) = \omega_0 = \const$, then Cayley-sine and Cayley-cosine functions satisfy the above dynamic equation. Therefore, this $\T$-analogue of the harmonic oscillator has \emph{bounded} solutions.

\section{Pad\'e-exponential functions on $\T$} The Pad\'e approximation consists in approximation by rational functions of prescribed order (see, e.g., \cite{An}).  Pad\'e approximants  $R_{j,k}$ are rational functions
\[
R_{j,k} (x) = \frac{P_{j} (x)}{Q_k (x)} \ ,
\]
which agree with $e^x$ (at $x=0$) to the highest possible order. Thus, e.g., $R_{j,k} (0) = 1$. 

\begin{definition}
The Pad\'e-exponential function $E^\alpha_{j,k}(t, t_0)$ satisfies the dynamic system defined at right dense points by $x^\Delta = \alpha x$ and at right scattered points by
\[
  x^\sigma = R_{jk} (\alpha \mu) \, x \ . 
\]
\end{definition}

\no Known exponential functions can be considered as particular Pad\'e exponentials: 
\[ \ba{lll}
 \dis E^\alpha_{1,0} (t, t_0) = e_\alpha (t, t_0) \quad & \emph{delta exponential function} \ , & \dis x^\sigma = (1 + \alpha \mu ) x  \ , \\[3ex]
 \dis E^\alpha_{0,1} (t, t_0) = {\hat e}_\alpha (t, t_0) \quad  & \emph{nabla exponential function}   \ ,  & \dis x^\sigma = \frac{1}{1 - \alpha \mu}  \ , \\[3ex]
 \dis E^\alpha_{1,1} (t, t_0) = E_\alpha (t, t_0) \quad & \emph{Cayley-exponential function}  \ ,  & \dis x^\sigma = \frac{1 + \frac{1}{2} \alpha \mu}{1 - \frac{1}{2} \alpha \mu}  \ . 
\ea \]
Another special case, $x = E^\alpha_{2,2} (t, t_0)$, satisfies equations: 
\[
x^\sigma = \frac{1 + \frac{1}{2} \alpha \mu + \frac{1}{12} (\alpha \mu)^2}{1 - \frac{1}{2} \alpha \mu + \frac{1}{12} (\alpha \mu)^2} \, x  , \quad 
 \dis x^\Delta = \frac{\alpha }{ 1 + \frac{1}{12} (\alpha \mu)^2 } \, \av{x} \equiv  
 \frac{\alpha  }{ 1 - \frac{1}{2} \alpha \mu + \frac{1}{12} (\alpha \mu)^2 } \, x  . 
\]
Symmetric Pad\'e-exponentials ($ E^\alpha_{k, k} (t^\sigma, t_0)$)  map the imaginary axis into the unit circle and, therefore, they generate trigonometric functions with ``good'' properties.

\section{ Exact analogues of elementary/special functions on $\T$ }
\label{sec-exact}

Given \ $f : \R \rightarrow \C$, we define its \emph{exact analogue} \  $\tilde f : \T \rightarrow \C$ \ as \  $\tilde f := f|_{\T}$, \ i.e.,  
\[
    \tilde f (t) := f (t) \quad ({\rm for}\ t \in \T) \ .
\]
The path   $f \rightarrow \tilde f$ is obvious and unique, but the inverse way (to find  $f$  corresponding to a given $\tilde f$) is, in general, neither obvious nor unique. However, in some particular cases a natural correspondence $\tilde f \rightarrow f$ may exist. The constant function $\tilde f$ is a typical example \cite{Ci-trig}. The corresponding $f$ is, obviously, constant as well.

\begin{definition} We assume $\alpha = \const \in \C$. 
 The \emph{exact exponential function} on  $\T$  is defined as 
\ $E^{ex}_\alpha (t, t_0) := e^{\alpha (t-t_0)}$. 
\end{definition}

Similarly we can define exact trigonometric and exact hyperbolic functions. Exact trigonometric functions on $\T$ coincide with trigonometric functions on $\T$ introduced by Hilger \cite{Hi-spec}, see also section~\ref{sec-hyptrig}. 

\begin{theorem} 
The exact exponential function $E^{ex}_\alpha (t, t_0)$ satisifies
\[
  x^\Delta (t)  = \alpha \ \psi_{\alpha} (t)  \  \av{x (t)} \ , 
\qquad x (t_0) = 1 \ ,
\]
where $\psi_\alpha (t) = 1$ at right-dense points and 
\[
  \psi_\alpha (t) = \frac{2}{\alpha \mu (t) }{\tanh \frac{\alpha \mu (t)}{2} }
\]
at right-scattered points.  
\end{theorem}

Exact discretizartions of ODEs with constant coefficients were first considerd by Potts \cite{Po}, see also \cite{Ci-oscyl,CR-ade,Mic}.

\section{A modification of the $q$-calculus} 

The C-exponential function suggests analogical modifications in the $q$-calculus. A survey of classical results can be found in the textbook \cite{KC}. It is curious that among several definitions of $q$-trigonometric functions, used throughout all the previous century, none satisfied the Pythagorean identity.  Below we present our recent results, see \cite{Ci-q}, enjoying this property. 

New $q$-exponential function can be defined in two equivalent ways, either as an infinite product or an infinite series:
\[
{\mathcal E}_q^x = \prod_{k=0}^\infty \frac{1 + q^k (1-q) \frac{x}{2} }{1 - q^k (1-q) \frac{x}{2}} \ , \qquad \quad {\mathcal E}_q^x = \sum_{n=0}^\infty  \frac{x^n}{\{n\}!} \ , 
\]  
where 
\[ 
 \{n\}! = \{1\} \{2\} \ldots \{n\}  \ , \qquad  \quad   \{k\} :=  \frac{1+q + \ldots + q^{k-1}}{ \frac{1}{2} (1 + q^{k-1}) } \ .
\] 
Our new $q$-exponential function can be directly expressed by  standard $q$-exponential functions:  
 ${\mathcal E}_q^x := e_q^{\frac{x}{2} } \, E_q^{\frac{x}{2}}  $  (for definitions and more details on  
 $e_q^x$ and $E_q^x$ see \cite{Ci-q,KC}).

New $q$-trigonometric functions,  motivated by the Cayley transformation, are defined by 
\[  \ba{l}  \dis
{\mathcal S}in_q x = \frac{{\mathcal E}_q^{i x} - {\mathcal E}_q^{-ix}}{2 i} \ , \quad {\mathcal C}os_q x = \frac{{\mathcal E}_q^{i x} + {\mathcal E}_q^{-ix}}{2 } \ ,
\ea \]
and have the following properties: 
\[  
 {\mathcal C}os_q^2 x + {\mathcal S}in_q^2 x = 1 \ , \qquad
    D_q {\mathcal S}in_q x = \av{ {\mathcal C}os_q x } \ , \qquad 
  D_q {\mathcal C}os_q x = - \av{ {\mathcal S}in_q x } \ ,
\]
where $q$-derivative is defined by $\dis D_q f (x) = \frac{f (qx) - f (q)}{q x - x}$ \ and \ 
$\dis \av{ f (x)} = \frac{f (x) + f (q x) }{2}$.

\section{Dynamic systems on Lie groups}

Two natural generalizations of the Cayley transformation are well known: 
 \begin{itemize}   
\item Lie algebra ${\mathfrak g}$ $\rightarrow$ (``quadratic'') Lie group $G$,   
\item anti-Hermitean operators  $\rightarrow$ unitary operators.
\end{itemize}
Here we focus on matrix Lie groups. It is well known that for quadratic Lie groups (including all orthogonal, unitary and symplectic groups) we have
\[
A \in {\mathfrak g} \quad  \Longrightarrow \quad (I - A)^{-1} (I + A) \in G \ . 
\]
Therefore, a natural $\T$-analogue of \ $\dis \frac{d}{d t} \Phi = A \Phi$ \ (here  $A \in {\mathfrak g }, \ \Phi  \in G $) is the dynamic equation \  $\Phi^\Delta = A \, \av{\Phi}$. The corresponding evolution at right-scattered points is expressed by the Cayley   transformation:
\[ 
 \Phi^\sigma = \frac{I + \frac{1}{2} \mu A}{I-\frac{1}{2} \mu A} \Phi  \ .
\]
A different approach to dynamic systems on G=$SU(2)$ can be found in \cite{Ci-pseudo}.

\section{Conclusions and future directions} 

The most important message of this paper consists in showing that one can define many non-equivalent exponential functions on time scales and many non-equivalent $\T$-analogues of ordinary differenatial equations. 
Differential equations have no unique `natural' time scales analogues. It is worthwhile to consider different numerical schemes  in this  context, compare \cite{Poe}.  

Our results suggest numerous developments. We name only few of them: 
dynamic systems preserving integrals of motion and Lyapunov functions, 
new developments in the $q$-calculus (e.g., modifications of $q$-gamma function and of the Jackson integral), modifications of $q$-Laplace and $q$-Fourier transformations, and locally exact $\T$-analogues of elementary functions, compare \cite{CR-grad}.

\subsection*{Acknowledgments} I would like to thank Stefan Hilger and Martin Bohner for inviting me to  the Special Session ``Differential, Difference and Dynamic Equations'' of 8th AIMS  International Conference on Dynamical Systems, Differential Equations and Applications, and for stimulating comments and discussions.

\end{document}